\documentclass[11pt,a4paper]{amsart}
\usepackage[utf8x]{inputenc}
\usepackage[english]{babel}
\usepackage{amsmath}
\usepackage{amsfonts}
\usepackage{amssymb}
\usepackage{graphicx}
\usepackage{xcolor}
\usepackage{lmodern}
\usepackage{mathrsfs}
\usepackage{graphicx,amssymb, verbatim,enumerate,amsmath,mathtools}
\usepackage{esint}
\usepackage{amssymb, amsthm}
\usepackage{enumerate, mathrsfs, hyphenat, graphicx}
\usepackage[font=sf, labelfont=sf]{caption}
\usepackage{thmtools, thm-restate}
\usepackage{thm-patch, thm-kv, thm-autoref}
\usepackage[colorlinks=true, urlcolor=blue, linkcolor=purple, citecolor=green]{hyperref}
\usepackage[left=2cm,right=2cm,top=2cm,bottom=2cm]{geometry}
\theoremstyle{plain}
\declaretheorem[title=Theorem, parent=section]{sa}
\declaretheorem[title=Lemma,sibling=sa]{lem}
\declaretheorem[title=Corollary,sibling=sa]{cor}
\declaretheorem[title=Proposition,sibling=sa]{prop}
\newtheorem*{thm*}{Satz}

\theoremstyle{definition}
\declaretheorem[title=Definition,sibling=sa]{defi}
\declaretheorem[unnumbered,title=Remark]{rem}

\numberwithin{equation}{section}
\newcommand{\al}{\alpha}
\def\Mitaga1{{E_{\al,1}}}

\newcommand{\R}{  \mathbb{R}}

\newcommand{\fdabc}[2]{\prescript{\text{abc}}{#1} D_{t}^{#2}}
\newcommand{\Fdabc}[2]{\prescript{\text{abc}}{#1} D_{T-t}^{#2}}
\newcommand{\fdabr}[2]{\prescript{\text{abr}}{#1} D_{t}^{#2}}

\begin{document}
\title[Optimal control of diffusion equation with fractional time derivative]{Optimal control of diffusion equation with fractional time derivative with nonlocal and nonsingular Mittag-Leffler kernel}
\author[Djida]{J.D. Djida}
\author{G.M.~Mophou}
\author{I.~Area}
\begin{abstract}
In this paper, we consider a diffusion equation with fractional-time derivative with nonsingular Mittag-Leffler kernel in Hilbert spaces. Existence and uniqueness of solution are proved by means of a spectral argument. The existence of solution is obtained for all values of the fractional parameter $\alpha \in (0,1)$. Moreover, by applying control theory to the fractional diffusion problem we obtain an optimality system which has also a unique solution. 
\end{abstract}
\address[Djida]{Departamento de An{\'a}lise Matem\'{a}tica, Universidade de Santiago de Compostela, 15782 Santiago de Compostela, Spain and African Institute for Mathematical Sciences (AIMS), P.O. Box 608, Limbe Crystal Gardens, South West Region, Cameroon}
\email[Djida]{\tt jeandaniel.djida@aims-cameroon.org}

\address[Mophou]{African Institute for Mathematical Sciences (AIMS), P.O. Box 608, Limbe Crystal Gardens, South West Region, Cameroon and Laboratoire C.E.R.E.G.M.I.A., D\'{e}partement de Math\'{e}matiques et Informatique, Universit\'{e} des Antilles et de la Guyane, Campus Fouillole, 97159 Pointe-\`{a}-Pitre,
(FWI), Guadeloupe}
\email[Mophou]{\tt gisele.mophou@aims-cameroon.org}

\address[Area]{Departamento de Matem\'atica Aplicada II,
E.E. Aeron\'autica e do Espazo, Universidade de Vigo,
Campus As Lagoas s/n, 32004 Ourense, Spain.}
\email[Area]{\tt area@uvigo.es}

\thanks{
\noindent The first author grateful for the facilities provided by the German research Chairs and the Teacher Training Program of AIMS-Cameroon. The first author is also indebted to the AIMS-Cameroon 2017-2018 Tutor fellowship. The second author was supported by the Alexander von Humboldt foundation, under the programme financed by the BMBF entitled "German research Chairs". The work of the third author has been partially supported by the Agencia Estatal de Innovaci\'on (AEI) of Spain under grant MTM2016-75140-P, cofinanced by the European Community fund FEDER, and Xunta de Galicia, grants
GRC 2015-004 and R 2016/022.}
\maketitle

\pagestyle{headings}		
\markboth{}{Optimal control of diffusion equation with fractional time derivative}

\section{Introduction}
\noindent In some problems related with transport dynamics in complex systems, the anomalous diffusion appear in a natural way \cite{Metzler}. As a consequence, some classical differential equations have been studied in the context of fractional derivatives, such as diffusion equation, diffusion-advection or Fokker-Plank type \cite{Metzler}. In \cite{Schneider} diffusion and wave equations together with appropriate initial condition(s) have been rewritten as integro-differential equations with time derivatives replaced by convolution with $t^{\alpha-1}/\Gamma(\alpha)$, $\alpha = 1,2$, respectively. These equations have been solved numerically by using different approaches (see e.g. \cite{Lei,Hao,Wang} and references therein). In \cite{Mophou} the classical control theory has been applied to a fractional diffusion equation in a bounded domain, where the fractional time derivative was considered in a Riemann-Liouville sense.

\medskip
\noindent We are interested in the analysis of the initial boundary value problem for a fractional time diffusion problem. Let $\Omega$ be an open and bounded subset of $\R^{d} (d\geq 1)$, with sufficiently smooth boundary $\partial \Omega $. Given $ \alpha \in (0,1)$, a well behave function $f \in L^{2}(Q)$, and an initial datum $y^{0} \in L^{2}(\Omega)$, we seek $y$ such that
\begin{equation}\label{Eq:problem}
\left\{
\begin{array}{ll}
\fdabc{0}{\alpha} y(x,t)-(\mathcal{A}y)(x,t) = f, & \text{for } (x,t)\in Q = \Omega \times (0,T),\\
y(x,0) = y^{0}(x), &\text{for } x \in \Omega,\\
y(x,t) = 0, & \text{on} ~\partial \Omega \times (0,T),
\end{array}
\right.
\end{equation}
where $\fdabc{a}{\alpha}$ for $\alpha \in (0,1)$ denotes the Atangana-Baleanu fractional derivative of order $\alpha$ in the sense of Caputo \cite{atanbaleanu} with respect to the time $t$, which is formally defined in \autoref{sec:prelim}.
\noindent In the cases of single and multiple fractional time Caputo derivatives, similar problems to \eqref{Eq:problem} have been already studied. For the Caputo fractional derivative there exists a large and rapidly growing number of publications (see e.g. \cite{Bazhlekova, Diethelm, Yamamoto, Mophou, Luchko,YuryLuchko} and references therein). Moreover, in \cite{Mophou} the Caputo fractional derivative was used to obtain an optimality system for the optimal control problem. \medskip

\noindent We would like to emphasize here that the existence and uniqueness of solutions for the fractional diffusion equation of type \eqref{Eq:problem} when the Caputo fractional-time derivative is considered, have been obtained for $\alpha \in (1/2,1)$ \cite{Yamamoto,Mophou,BaleanuJosephMophou}. In this work using the fractional-time derivative with the nonsingular Mittag-Laffler kernel, we have obtained the existence and uniqueness results for the fractional diffusion equation of type \eqref{Eq:problem} for all the values $\alpha \in (0,1)$. \medskip

\noindent Furthermore, as for as the optimal control of the problem \eqref{Eq:problem}, one can refer to the methods of the Lagrange multiplier technique for the Caputo and Riemann fractional time derivative presented by different authors (see e.g. \cite{Yamamoto, Mophou, BaleanuJosephMophou, Agrawal, Agrawal1, Delfim}  and references therein). We are interested in finding the control $\hat{u} \in L^{2}(Q)$ that minimizes the cost function
\[
\mathcal{J}(v) = \big\|y(v)-z_{d}\big \|^{2}_{L^{2}(\mathcal{Q})} + \frac{\mathcal{N}}{2}\|v\|^{2}_{L^{2}(Q)}, \quad z_{d} \in L^{2}(Q), ~~ \mathcal{N}>0,
\]
subject to the problem \eqref{Eq:problem} where the the Atangana-Baleanu fractional-time derivative is involved. In other words we are interested in finding the control parameter $\hat{u} \in L^{2}(Q)$ such that the functional
\[
\mathcal{J}(\hat{u}) = \inf_{v\in L^{2}(Q)} \mathcal{J}(v),
\]
holds true.\medskip

\noindent This manuscript is structured as follows: In section \ref{sec:prelim}  we collect notation, definitions, and preliminary results regarding the solution representation and we compute the integration by parts involving the Atangana-Baleanu fractional-time derivative. Consequently the weak formulation of the \eqref{Eq:problem} is obtained.  Section \ref{sec:existence}, is dedicated to the sketch of proof of the existence and uniqueness of solution of the weak formulation problem. More precisely we prove Theorem \ref{existence_uniqueness} as well as some corollaries. Moreover, section \ref{control} is devoted to an optimal control problem given by Theorem \ref{optimality_system}.

\section{Preliminaries}\label{sec:prelim}

\noindent In this section, we collect some notations on the functional spaces, useful properties on the Mittag-Leffler function and the fractional time derivatives, regularity results for the fractional diffusion \eqref{Eq:problem}, and basic estimates.\medskip

\subsection{Notations}\label{notations}
\noindent Let $L^{2}$ be the usual Hilbert space equipped with the scalar product $(\cdot,\cdot)$, and $H^{\ell}(\Omega)$, $H_{0}^{m}(\Omega)$ denote the
usual Sobolev spaces. In order to make things clear, we define the second order elliptic operator as
\begin{equation}
\label{second_order}
 \mathcal{A}y := -\Delta y, \text{ for}~ x \in \overline{\Omega}, \quad \mathcal{D}(A) = H^{2}(\Omega) \cap H_{0}^{1}(\Omega).
\end{equation}

\noindent Since $\mathcal{A}$ is an uniformly elliptic operator, the spectrum of $\mathcal{A}$ is composed of
eigenvalues $\big\{\lambda_{n}\big\}_{n=1}^{\infty}$ and its corresponding orthogonal eigenfunctions
$\big\{\varphi_{n}\big\}_{n=1}^{\infty}$ which belong to $H^{2}(\Omega) \cap H_{0}^{1}(\Omega)$, such that
$\mathcal{A}\varphi_{n} = \lambda_{n} \varphi_{n}$.

\subsection{Atangana-Baleanu fractional times derivatives: definitions and some properties}\label{prelim}

\noindent We first start by recalling the Mittag-Leffler function $E_{\alpha,\beta}(z)$ which will be used extensively throughout this work and is defined below
\[
E_{\alpha,\beta}(z) = \sum_{k=0}^\infty \frac{z^k}{\Gamma(k\alpha+\beta)},\quad z\in \mathbb{C},
\]
where $\Gamma(\cdot)$ denotes the Gamma function defined as
\[
\Gamma(z) = \int_0^\infty t^{z-1}e^{-t} dt , \quad \Re(z) >0.
\]

\noindent The Mittag-Leffler function is a two-parameter family of entire functions of $z$ of order $\al^{-1}$ \cite[pp. 42]{Kilbas}.
The exponential function is a particular case of the Mittag-Leffler function, namely $E_{1,1}(z) = e^z$, \cite{Kilbas,Samko}.
As presented in \cite{Jin}, two important members of this family are $E_{\alpha,1}(-\lambda t^\alpha)$ and $t^{\alpha-1}E_{\alpha,\alpha}(-\lambda t^\alpha)$, which occur in the solution operators for the initial value problem \eqref{Eq:problem}.\medskip

\noindent Next we also recall the generalized Mittag-Leffler function \cite{Kilbas, Samko} defined as
\begin{equation}\label{GIO}
\mathcal{E}^{\rho}_{\alpha,\beta}(z)=\sum_{k=0}^{\infty} \frac{(\rho)_{k}z^{k}}{\Gamma(\alpha k + \beta) k!},\qquad \forall t\in \mathbb{C},
\end{equation}
where $(\rho)_{k}$ is the  Pochhammer symbol defined as $(\rho)_{k} = \Gamma(\rho + k)/\Gamma(\rho) = \rho(\rho + 1)\cdots (\rho + k -1)$.
\begin{rem}
Once could notice that for $\rho =1$ we have that
\[
\mathcal{E}^{1}_{\alpha,\beta}(z)= E_{\alpha,\beta}(z).
\]
\end{rem}
\noindent Furthermore we recall the following Lemma from \cite{Kilbas}.
\begin{lem}\label{generalized}
Let $\alpha,\beta, \rho \in \mathbb{C}$ such that $\Re(\alpha) >0$ and $\Re(\beta)> 0$. Then we have that
\begin{equation}\label{first}
\left(\frac{d}{dz} \right)^{k} \mathcal{E}^{\rho}_{\alpha,\beta}(z) = (\rho)_{k}~ \mathcal{E}^{\rho+k}_{\alpha,\alpha k + \beta}(z), ~~ z \in \mathbb{C},~~ k \in \mathbb{N},
\end{equation}
and in the other hand,
\begin{equation}\label{ssecond}
\alpha \rho ~ \mathcal{E}^{\rho + 1}_{\alpha, \beta}(z) = (1 + \alpha \rho - \beta) \mathcal{E}^{\rho}_{\alpha,\beta}(z) + \mathcal{E}^{\rho}_{\alpha, \beta-1}(z), ~~ z \in \mathbb{C}.
\end{equation}
\end{lem}
\noindent Let us recall some useful definitions of fractional derivatives in the sense of Atangana-Baleanu \cite{atanbaleanu, Djida}.\medskip

\begin{defi}\label{def_atanbaleanu_C}
For a given function $u \in H^{1}(a,T)$, $T>a$, the Atangana-Baleanu fractional derivative of $u$ of order $\alpha$ in Caputo sense  with base
point $a$ is defined at a point $t \in (a,T)$ by
\begin{equation}\label{eq:abfdc}
\fdabc{a}{\alpha}~u(t)= \frac{B(\alpha)}{1-\alpha} \int_{a}^{t} u'(\tau) E_{\alpha} \left[ -\gamma(t-\tau)^{\alpha} \right] d\tau,
\end{equation}
where $\gamma = \alpha/(1-\alpha)$, $E_{\alpha}(\cdot)$ stands for the Mittag-Leffler function, and $B(\alpha)= (1-\alpha) +\alpha /\Gamma (\alpha )$.
\end{defi}

So from the definition on \cite{ThabetBaleanu} we recall the following definition
\begin{defi}\label{def_atanbaleanu_C_at_T}
For a given function $u \in H^{1}(a,T)$, $T>t>a$, the Atangana-Baleanu fractional derivative of $u$ of order $\alpha$ in Caputo sense  with base
point $T$ is defined at a point $t \in (a,T)$ by

\begin{equation}\label{Shiftatanbaleanu}
\fdabc{T}{\alpha}~g(t)= -\frac{B(\alpha)}{1-\alpha} \int_{t}^{T} g'(\tau) E_{\alpha} \left[ -\gamma(\tau-t)^{\alpha} \right] d\tau,
\end{equation}
\end{defi}

\begin{defi}\label{def_atanbaleanu_R}
\noindent Let $u \in H^{1}(a,T)$, $T>a$, $\alpha \in (0,1)$. The Atangana-Baleanu fractional derivative of order $\alpha$ of $u$ in Riemann-Liouville sense with base point $a$ is defined at a point $t \in (a,T)$ as
\begin{equation}\label{eq:fdab}
\fdabr{a}{\alpha} u(t)=\frac{B(\alpha)}{1-\alpha} \frac{d}{dt} \int_{a}^{t} u(\tau) E_{\alpha} \left[ -\gamma (t-\tau)^{\alpha}\right] d\tau.
\end{equation}
\end{defi}

For $\alpha = 1$ in \eqref{eq:abfdc}, we consider the usual classical derivative $\partial_t$.\medskip

\noindent Furthermore the Atangana-Baleanu fractional integral of order $\alpha$ with base point $a$ is defined as
\begin{equation}\label{eq:fi}
I_{t}^{\alpha} u(t) = \frac{1-\alpha}{B(\alpha)}u(t) + \frac{\alpha}{B(\alpha) \Gamma(\alpha)} \int_{a}^{t} u(\tau) (t-\tau)^{\alpha-1} d\tau.
\end{equation}
Notice that if $\alpha=0$ in \eqref{eq:fi} we recover the initial function, and if $\alpha=1$ in \eqref{eq:fi} we obtain the ordinary integral.

\noindent Next we recall the following lemma which has been stated and proved in \cite[Theorem $1.6$]{Podlubny}
\begin{lem}\label{lem:mlfbdd}
Let $\alpha \in (0,1)$ and $\beta\in\mathbb{R}$ be arbitrary, and $\mu$ is such that $\frac{\alpha\pi}{2}
<\mu<\min(\pi,\alpha\pi)$. Then there exists a real constant $C=C(\alpha,\beta,\mu)>0$ such that
\begin{equation}\label{M-L-bound}
|E_{\alpha,\beta}(z)|\leq \frac{C}{1+|z|},\quad\quad \mu\leq|\mathrm{arg}(z)|\leq \pi.
\end{equation}
\end{lem}

\begin{lem}\label{shift}
 Let set $\eta : [0,T] \to \mathbb{R}$, then for all $\alpha \in (0,1)$ we the equivalence relation hold true
\begin{equation}\label{shift_operator}
 \fdabc{T}{\alpha}~\eta(T-t) =  \fdabc{0}{\alpha}~\eta(t).
\end{equation}
\end{lem}
\begin{proof}
 The well-posedness of \eqref{shift_operator} follows from a change of variables. Let us set $\eta : [0,T] \to \mathbb{R}$, and let us define by
 $ \widetilde{\eta}(t):= \eta(T-t)$. Notice that
 \[
  \widetilde{\eta}'(t) = -\eta'(T-t).
 \]

\noindent Next by making change of variables $t \to T-t$ in \eqref{Shiftatanbaleanu}, we get
\[
\Fdabc{0}{\alpha}~\eta(T-t) = \fdabc{T}{\alpha}~\widetilde{\eta}(t)
 = -\frac{B(\alpha)}{1-\alpha} \int_{T-t}^{T} \eta'(\tau) E_{\alpha} \left[ -\gamma(\tau-(T-t))^{\alpha} \right] d\tau,
\]
setting $\omega = T-\tau$, we get that
\[
\begin{aligned}
 &\fdabc{T}{\alpha}~\widetilde{\eta}(t) = \frac{B(\alpha)}{1-\alpha} \int_{T-t}^{T} \eta'(T-\omega) E_{\alpha} \left[ -\gamma(t-\omega)^{\alpha} \right] d\omega \\
 &= \frac{B(\alpha)}{1-\alpha} \int_{0}^{t} \eta'(T-\omega) E_{\alpha} \left[ -\gamma(t-\omega)^{\alpha} \right] d\omega
\fdabc{T}{\alpha}~\widetilde{\eta}(t) \\
&=- \frac{B(\alpha)}{1-\alpha} \int_{0}^{t} E_{\alpha} \left[ -\gamma(t-\omega)^{\alpha} \right]\widetilde{\eta}'(\omega ) d\omega,
\end{aligned}
\]
This means that
\[
\fdabc{T}{\alpha}~\widetilde{\eta}(t) = -\fdabc{0}{\alpha}~\widetilde{\eta}(t)
\]
\end{proof}
\noindent As a consequence, the backwards in time with the fractional-time derivative with nonsingular Mittag-Leffler kernel at the based point
$T$ is equivalently written as a forward in time operator with the fractional-time derivative with nonsingular Mittag-Leffler kernel
$-\fdabc{0}{\alpha}$.

\subsection{Solution representation }\label{sol}
\noindent Under the functional spaces introduced in \autoref{notations}, we can associate with $u(x,t)$ a map $u(\cdot):(0,T) \to L^{2}(\Omega)$ by
\[
u(t)(x) = u(x,t), \quad 0<t<T, ~~ x \in \Omega.
\]

Let us now multiply \eqref{Eq:problem} by a function $\vartheta \in H^{1}_{0}(\Omega)$ and by integration by parts on $\Omega$, we get that
\[
\int_{\Omega} \fdabc{0}{\alpha} u(x,t)\vartheta (x) dx - \int_{\Omega} \Delta u(x,t) \vartheta (x)dx = \int_{\Omega} f(x,t)\vartheta (x)dx.
\]
By Green's formula, we have
\begin{equation}\label{Eq:repre}
\int_{\Omega} \fdabc{0}{\alpha} u(x,t)\vartheta (x) dx + \int_{\Omega} \nabla u(x,t) \nabla \vartheta (x)dx = \int_{\Omega} f(x,t)\vartheta (x)dx.
\end{equation}
Now we set for all \[
\varphi, \psi \in L^{2}(\Omega),~~ (\varphi,\psi)_{L^{2}(\Omega)} = \int_{\Omega}\varphi(x)\psi(x) dx.
\]
\medskip
\noindent The scalar product in $L^{2}(\Omega)$, with endowed norm $\|\varphi\|_{L^{2}(\Omega)}$. In the same spirit, we define the bilinear form in $H_{0}^{1}(\Omega)$ as
\begin{equation}\label{bilinear_form}
a(\varphi,\psi) = \int_{\Omega} \nabla \varphi(x) \nabla \psi (x)dx, \quad \forall \varphi, \psi \in H^{1}_{0}(\Omega),
\end{equation}
so that we have $\|\varphi\|^{2}_{H^{1}_{0}(\Omega)} = a(\varphi,\psi)$.\medskip

\noindent Under the assumption that the operator $\mathcal{A}$ is an uniformly elliptic operator, the spectrum of $\mathcal{A}$ is composed of eigenvalues $\{\lambda_{k}\}_{k=1}^{\infty}$, and its corresponding orthogonal eigenfunctions $\{\lambda_{w_{k}}\}_{k=1}^{\infty}$ which belongs to $\mathcal{D}(A)$, such that
\[
a(w_{k},p) = \lambda_{k}(w_{k},r)_{L^{2}(\Omega)}, \quad \forall r \in ~H_{0}^{1}(\Omega).
\]
Hence we have in this setting
\begin{equation}\label{Eq:repre2}
\|\varphi\|^{2}_{H^{1}_{0}(\Omega)} = \sum_{k=0}^{\infty}\lambda_{i}(\varphi,w_{i})^{2}_{L^{2}(\Omega)}, \quad \forall \varphi \in ~H_{0}^{1}(\Omega).
\end{equation}

\noindent As an immediate consequence, equation \eqref{Eq:repre2} can be written as follows
\[
(\fdabc{0}{\alpha} u(x,t),\vartheta )_{L^{2}(\Omega)} + a(u(t),~ \vartheta)_{L^{2}(\Omega)}  = (f(t),~\vartheta)_{L^{2}(\Omega)}.
\]
In this way, the problem \eqref{Eq:problem} takes the representation as

\begin{equation}\label{Eq:newproblem}
\left\{
\begin{array}{ll}
\big(\fdabc{0}{\alpha} u(x,t),\vartheta \big)_{L^{2}(\Omega)} + a\big(u(t),~ \vartheta\big)_{L^{2}(\Omega)}  = \big(f(t),~\vartheta \big)_{L^{2}(\Omega)} & \text{in} ~~\Omega,~~ \forall ~\vartheta \in H_{0}^{1}(\Omega),\\
u(x,0) = u^{0} &\text{for } x \in \Omega,\\
u(t) = 0 & \text{on} ~\partial \Omega.
\end{array}
\right.
\end{equation}


\noindent To discuss the existence of uniqueness results of \eqref{Eq:problem} or \eqref{Eq:newproblem}, we shall need some  informations over the various time. We are interested in the solution of the fractional differential equation involving the Atangana-Baleanu fractional derivative, in the form of \eqref{Eq:newproblem}. To this end, we propose to solve the fractional integro-differential equation involving the Atangana-Baleanu fractional derivative in the form
\begin{equation}\label{Eq:fdfeq}
\left\{
\begin{array}{ll}
\fdabc{0}{\alpha} y_{i}(t) + \lambda_{i}y_{i}(t) = f_{i}(t), \qquad t\in (0,T),\\[3mm]
y_{i}(0) = y_{i}^{0}.
\end{array}
\right.
\end{equation}
\noindent This type of fractional differential equation involving the fractional-time derivative with nonsingular Mittag-Leffler kernel has been solved in \cite[Proposition $3.3$]{DjidaIvan_Abdon_regularity} by the means of Laplace transform. By following the same procedure, we state the the next proposition
\begin{prop}\label{proposition_solution}
\noindent Let $f \in L^{2}(\mathcal{Q})$, such that the Atangana-Baleanu fractional derivative exists. Then, the solution of differential equation \eqref{Eq:fdfeq}, is given by
\begin{equation}\label{Eq:sol2}
\begin{aligned}
y_{i}(t)  = \zeta_{i} E_{\alpha}\big[-\gamma_{i} t^{\alpha}\big]y_{i}^{0}
+ \frac{(1-\alpha) \zeta_{i}}{B(\alpha)}f_{i}(t) + K_{i} \int_{0}^{t}(t-s)^{\alpha} E_{\alpha,\alpha}\big[-\gamma_{i}(t-s)^{\alpha}\big]f_{i}(s)ds,
\end{aligned}
\end{equation}
with
\begin{equation}\label{contants}
\gamma_{i} = \frac{\alpha \lambda_{i}}{\big( B(\alpha) + (1-\alpha)\lambda_{i} \big)},\quad \zeta_{i} = \frac{B(\alpha)}{\big( B(\alpha) + (1-\alpha)\lambda_{i} \big)},\quad \text{and}\quad K_{i} = \bigg(\frac{\alpha \zeta_{i}}{B(\alpha)\Gamma(\alpha)} + \frac{(1-\alpha) \gamma_{i}\zeta_{i}}{B(\alpha)} \bigg).
\end{equation}
\end{prop}

\subsection{Weak formulation of the problem}\label{weak_formulation}

\noindent Next we state the following proposition which gives the weak formulation of the problem \eqref{Eq:problem}, that will be fundamental in our analysis.

\begin{prop}\label{weak_formulation_lem}
\noindent Let $\phi, y \in \mathcal{C}^{\infty}(\bar{Q})$. Then, we have
\begin{equation}\label{Eq:weak}
\begin{aligned}
&\int_{\Omega}\int_{0}^{T} \bigg(\fdabc{0}{\alpha}~y(x,t) - \Delta y(x,t) \bigg) \phi(x,t) dx~dt  =
\int_{0}^{T} \int_{\partial \Omega} y \frac{\partial \phi}{\partial \sigma} \phi d\sigma~dt - \int_{0}^{T} \int_{\partial \Omega}
\phi \frac{\partial y}{\partial \sigma} d\sigma~dt \\
& + \int_{\Omega}\int_{0}^{T}y(x,t) \bigg(-\fdabc{T}{\alpha} ~\phi(x,t) - \Delta  \phi(x,t) \bigg)dt~dx
-\frac{B(\alpha)}{1-\alpha}\int_{\Omega}\int_{0}^{T}  y(x,0)E_{\alpha,\alpha}[-\gamma t^{\alpha}]\phi(x,t) dt~dx\\
& + \frac{B(\alpha)}{1-\alpha}\int_{\Omega} \phi(x,T)\int_{0}^{T} y(x,t) E_{\alpha,\alpha}\left[-\gamma(T-t)^{\alpha} \right] dt~dx.
\end{aligned}
\end{equation}
\end{prop}
\begin{proof}
For a given $\phi, y \in \mathcal{C}^{\infty}(\bar{Q})$, we have that
\[
\begin{aligned}
&\int_{\Omega}\int_{0}^{T} \bigg(\fdabc{0}{\alpha}~y(x,t) - \Delta y(x,t) \bigg) \phi(x,t) dx~dt  \\
& = \int_{\Omega}\int_{0}^{T} \fdabc{0}{\alpha}~y(x,t)\phi(x,t)dx~dt - \int_{\Omega}\int_{0}^{T} \Delta y(x,t)\phi(x,t) dx~dt  \\
&  = (A) + (B)
\end{aligned}
\]
Let us compute each part.\medskip

\noindent $\bullet$ {\sc computation of $(A)$: (Integration by parts)}\medskip

\noindent Using the definition of the Atangana-Baleanu fractional derivative \eqref{eq:abfdc} at the based point $a=0$ we have the following
\begin{equation}\label{computation_I_a}
\int_{0}^{T} \fdabc{0}{\alpha}~y(x,t)\phi(x,t)dt = \int_{0}^{T} \left[\frac{B(\alpha)}{1-\alpha} \int_{0}^{t} y'(\tau) E_{\alpha} \left[-\gamma(t-\tau)^{\alpha} \right] d\tau \right] \phi(t)dt.
\end{equation}
We define by
\[
A_{1} := \int_{0}^{t} y'(\tau) E_{\alpha} \left[-\gamma(t-\tau)^{\alpha} \right] d\tau
\]
Then by computing $A_{1}$ by integration by parts, we get that
\[
A_{1} = y(t) - y(0)E_{\alpha}[-\alpha t^{\alpha}] - \int_{0}^{t}\gamma(t-\tau)^{\alpha-1}E_{\alpha,\alpha} \left[-\gamma(t-\tau)^{\alpha} \right]y(\tau)d\tau.
\]
Replacing $A_{1}$ by its value in \eqref{computation_I_a} we get that
\begin{equation}\label{computation_I_ab}
\begin{aligned}
\int_{0}^{T} \fdabc{0}{\alpha}~y(x,t)\phi(x,t)dt &= \frac{B(\alpha)}{1-\alpha} \int_{0}^{T} y(x,t) \phi(x,t) dt - \frac{B(\alpha)}{1-\alpha} \int_{0}^{T} y(x,0)E_{\alpha}[-\alpha t^{\alpha}]\phi(x,t)dt \\
&- \frac{B(\alpha)}{1-\alpha} \int_{0}^{T}  \left[ \int_{0}^{t} \gamma (t-\tau)^{\alpha-1}E_{\alpha,\alpha} \left[-\gamma(t-\tau)^{\alpha} \right]y(x,\tau)d\tau \right]\phi(x,t)dt\\
&:= A_{2} + A_{3} + A_{4},
\end{aligned}
\end{equation}
with 
\[
\begin{aligned}
& A_{2}:= \frac{B(\alpha)}{1-\alpha} \int_{0}^{T} y(x,t) \phi(x,t) dt; \quad A_{3}:=\frac{B(\alpha)}{1-\alpha} \int_{0}^{T} y(x,0)E_{\alpha}[-\alpha t^{\alpha}]\phi(x,t)dt \\
&A_{4}:= \frac{B(\alpha)}{1-\alpha} \int_{0}^{T}  \left[ \int_{0}^{t} \gamma (t-\tau)^{\alpha-1}E_{\alpha,\alpha} \left[-\gamma(t-\tau)^{\alpha} \right]y(x,\tau)d\tau \right]\phi(x,t)dt.
\end{aligned}
\]

Let us compute $A_{4}$
\[
\begin{aligned}
A_{4} &:= -\frac{B(\alpha)}{1-\alpha} \int_{0}^{T}  \left[ \int_{0}^{t} \gamma (t-\tau)^{\alpha-1}E_{\alpha,\alpha} \left[-\gamma(t-\tau)^{\alpha} \right]y(x,\tau)d\tau \right]\phi(x,t)dt \\
&= \frac{B(\alpha)}{1-\alpha} \int_{0}^{T}  y(x,\tau)\left[ -\int_{\tau}^{T} \gamma (t-\tau)^{\alpha-1}E_{\alpha,\alpha} \left[-\gamma(t-\tau)^{\alpha} \right] \phi(x,t) dt \right]d\tau\\
&= \frac{B(\alpha)}{1-\alpha} \phi(x,T)\int_{0}^{T} y(x,t) E_{\alpha,\alpha} \left[-\gamma(T-t)^{\alpha} \right]dt -\frac{B(\alpha)}{1-\alpha} \int_{0}^{T} y(x,t) \phi(x,t) dt \\
&+ \int_{0}^{T} y(x,t) \left[ -\frac{B(\alpha)}{1-\alpha} \int_{t}^{T} E_{\alpha,\alpha} \left[-\gamma(t-\tau)^{\alpha} \right]\phi'(x,\tau) d\tau \right]dt \\
&= \frac{B(\alpha)}{1-\alpha} \phi(x,T)\int_{0}^{T} y(x,t) E_{\alpha,\alpha} \left[-\gamma(T-t)^{\alpha} \right]dt -\frac{B(\alpha)}{1-\alpha} \int_{0}^{T} y(x,t) \phi(x,t) dt \\
&+\int_{0}^{T}y(x,t) \fdabc{T}{\alpha}~y(x,t)\phi(x,t)dt
\end{aligned}
\]
Hence replacing $A_{4}$ by its value in \eqref{computation_I_ab} we finally have that
\begin{equation}\label{computation_I}
\begin{aligned}
&\int_{0}^{T} \fdabc{0}{\alpha}~y(x,t)\phi(x,t)dt = -\int_{0}^{T}y(x,t) \fdabc{T}{\alpha}~y(x,t)\phi(x,t)dt \\
&+ \frac{B(\alpha)}{1-\alpha} \phi(x,T)\int_{0}^{T} y(x,t) E_{\alpha,\alpha} \left[-\gamma(T-t)^{\alpha} \right]dt - \frac{B(\alpha)}{1-\alpha}  y(x,0)\int_{0}^{T}E_{\alpha}[-\alpha t^{\alpha}]\phi(x,t)dt
\end{aligned}
\end{equation}

\noindent $\bullet$ {\sc computation of $(B)$}\medskip

\noindent By a simple integration by parts we have that
\begin{equation}\label{computation_II}
\begin{aligned}
(B)&:=-\int_{\Omega}\int_{0}^{T}\phi(x,t)\Delta u(x,t)  dx~dt  = \int_{0}^{T} \int_{\partial \Omega} u \frac{\partial \phi}{\partial \sigma} \phi d\sigma~dt - \int_{0}^{T} \int_{\partial \Omega}  \phi \frac{\partial u}{\partial \sigma} d\sigma~dt\\
& - \int_{\Omega}\int_{0}^{T}u(x,t) \Delta  \phi(x,t)dt~dx.
\end{aligned}
\end{equation}

Hence adding \eqref{computation_I} to \eqref{computation_II}, we obtain the desired result \eqref{Eq:weak}.

\end{proof}
\noindent Now using the conditions that $\phi(x,T)=0$, we have the following Corollary.\medskip
\begin{cor}\label{corweak}
\noindent Let $y$ be the solution of \eqref{Eq:problem}. Then, for any $\phi, y \in  \mathcal{C}^{\infty}(\bar{Q})$ such that $\phi(x,T) = 0$ in $\Omega $, we have that,
\begin{equation}\label{Eq:weak_secong_version}
\begin{aligned}
&\int_{\Omega}\int_{0}^{T} \bigg(\fdabc{0}{\alpha}~y(x,t) - \Delta y(x,t) \bigg) \phi(x,t) dx~dt  =
\int_{0}^{T} \int_{\partial \Omega} y \frac{\partial \phi}{\partial \sigma} \phi d\sigma~dt - \int_{0}^{T} \int_{\partial \Omega}
\phi \frac{\partial y}{\partial \sigma} d\sigma~dt \\
& + \int_{\Omega}\int_{0}^{T}y(x,t) \bigg(-\fdabc{T}{\alpha} ~\phi(x,t) - \Delta  \phi(x,t) \bigg)dt~dx
 -\frac{B(\alpha)}{1-\alpha}\int_{\Omega}\int_{0}^{T} y(x,0) E_{\alpha,\alpha}[-\gamma t^{\alpha}]\phi(x,t) dt~dx.
\end{aligned}
\end{equation}
\end{cor}

\section{Existence and uniqueness of solution}\label{sec:existence}
\noindent Let us consider the following problem, for $\alpha \in (0,1)$. Assume $f \in L^{2}(Q)\big)$ and $y^{0} \in L^{2}(\Omega)$. We want to find $y \in L^{2}\big((0,T), H^{1}_{0}(\Omega)\big)$, $y(x,0) \in \mathcal{C}\big([0,T];H^{1}_{0}(\Omega)\big)$, such that $\forall ~~\vartheta \in H^{1}_{0}(\Omega)$, we have

\begin{equation}\label{Eq:newproblem1}
\left\{
\begin{array}{ll}
\big(\fdabc{0}{\alpha} y(t),\vartheta \big)_{L^{2}(\Omega)} + a\big(y(t),~ \vartheta\big)_{L^{2}(\Omega)}  = \big(f(t),~\vartheta \big)_{L^{2}(\Omega)} & \text{for all} ~~t \in (0,T),\\[3mm]
y(x,0) = y^{0} &\text{for } x \in \Omega,
\end{array}
\right.
\end{equation}
\medskip

\begin{sa}\label{existence_uniqueness}
\noindent Let $\alpha \in (0,1)$. Assume $f \in L^{2}(Q)\big)$, $y^{0} \in L^{2}(\Omega)$ and the bilinear form $a(\cdot,\cdot)$ defined as in \eqref{bilinear_form}. Then the problem \eqref{Eq:problem} or \eqref{Eq:newproblem1} has a unique solution $y \in \left(L^{2}\big((0,T); H_{0}^{1}(\Omega)\big) \right) \cap \left( \mathcal{C}\left([0,T];L^{2}(\Omega) \right) \right) $ given by
\[
\begin{aligned}
y(x,t)  = \sum_{i=1}^{+\infty}\bigg\{ \zeta_{i} E_{\alpha}\big[-\gamma_{i} t^{\alpha}\big]y_{i}^{0}
+ \frac{(1-\alpha) \zeta_{i}}{B(\alpha)}f_{i}(t) + K_{i} \int_{0}^{t}(t-s)^{\alpha-1} E_{\alpha,\alpha}\big[-\gamma_{i}(t-s)^{\alpha}\big]f_{i}(s)ds\bigg\}w_{i},
\end{aligned}
\]
where the constants $\gamma_{i}$ and $\zeta_{i}$ are defined as in \eqref{contants}. Moreover $y$ satisfies the bounded conditions
\begin{equation}\label{aaa}
\big\|y\big\|_{L^{2}\big((0,T);H_{0}^{1}(\Omega)\big)} \leq \Lambda_{1} \bigg( \|y^{0}\|_{H_{0}^{1}(\Omega)\big)} + \|f\|_{L^{2}(Q)}\bigg),
\end{equation}
where $\Lambda_{1} = \max\left(\mathcal{C}_{1}(\alpha,\lambda_{1},T),~ \mathcal{C}_{2}(\alpha,\lambda_{1},T)\right)$, with
\begin{equation}\label{ccc}
\mathcal{C}_{1}(\alpha,\lambda_{1},T) =  \frac{C B(\alpha)}{(1-\alpha)}\sqrt{\frac{6T}{\lambda_{1}}}, \qquad \mathcal{C}_{2}(\alpha,\lambda_{1},T) = \sqrt{\frac{6B^{2}(\alpha)}{(1-\alpha)^{2}\lambda_{1} } + \frac{12C^{2}T^{2}}{\Gamma^{2}(\alpha)} \left(\frac{\Gamma^{2}(\alpha)+ 1}{\lambda_{1}} \right)},
\end{equation}
and
\begin{equation}\label{aaaa}
\big\|y \big\|_{L^{2}(\Omega)} \leq \Lambda_{2} \left( \|y^{0}\|_{L^{2}(\Omega)} + \|f\|_{L^{2}(Q)}\right),
\end{equation}
where $\Lambda_{2} = \max\left(\mathcal{C}_{3}(\alpha,\lambda_{1},T),~ \mathcal{C}_{4}(\alpha,\lambda_{1},T)\right)$, with
\[
\mathcal{C}_{3}(\alpha,\lambda_{1},T) =  \frac{C B(\alpha)}{\lambda_{1}(1-\alpha)}\sqrt{6T}, \qquad \mathcal{C}_{4}(\alpha,\lambda_{1},T) = \sqrt{\frac{6}{\lambda_{1}^{2}} + \frac{12C^{2}T^{2}}{\Gamma^{2}(\alpha)} \left(\frac{\Gamma^{2}(\alpha)+ 1}{\lambda_{1}^{2}} \right)},
\]
providing that $y^{0}\in L^{2}(\Omega)$.
\end{sa}

\begin{proof}\label{frist_theorem}
\noindent We first assume that the series defines in \eqref{Eq:repre2} converges uniformly. Let $\vartheta = w_{i}$ in \eqref{Eq:newproblem1} then we have that
\begin{equation}\label{Eq:problem3}
\left\{
\begin{array}{ll}
\fdabc{0}{\alpha} y_{i}(t) + \lambda_{i}y_{i}(t) = f_{i}(t), \qquad t\in (0,T),\\[3mm]
y_{i}(0) = y_{i}^{0},
\end{array}
\right.
\end{equation}
where we have used
\[
a(y(t),w_{i}) = \lambda_{i}(y(t),w_{i})_{L^{2}(\Omega)} = \lambda_{i}y_{i}.
\]

\noindent We then recognize the system given by \eqref{Eq:problem3} which has been solved in subsection \ref{sol}. So the solution of the fractional differential equation takes the form
\begin{equation}\label{Eq:solution2}
\begin{aligned}
y_{i}(t)  = \zeta_{i} E_{\alpha}\big[-\gamma_{i} t^{\alpha}\big]y_{i}^{0}
+ \frac{(1-\alpha) \zeta_{i}}{B(\alpha)}f_{i}(t) + K_{i} \int_{0}^{t}(t-s)^{\alpha-1} E_{\alpha,\alpha}\big[-\gamma_{i}(t-s)^{\alpha}\big]f_{i}(s)ds.
\end{aligned}
\end{equation}
\noindent Notice that the $\lambda_{i}$ are contained in $\zeta_{i} $, $\gamma_{i}$ and $K_{i}$ previously defined in subsection \ref{sol}. \medskip

\noindent Now the next step is to look at the uniqueness of the solution. We will make it in two steps. We shall start by showing that the solution belongs to the space $L^{2}((0,T); H^{1}_{0}(\Omega))$ and in the second step, we will show that the solution also belong to $\mathcal{C}([0,T]; L^{2}(\Omega))$, so that $y $ belongs to $L^{2}((0,T); H^{1}_{0}(\Omega)) \cap \mathcal{C}([0,T]; L^{2}(\Omega)) $. \medskip

\noindent Let defined by $\mathcal{V}_{m}$ a subspace of $H_{0}^{1}(\Omega)$ generated by $w_{1},w_{2},\dots,w_{m}$. We want to find $y_{m}: t\in (0,T] \to y_{m}(t) \in \mathcal{V}_{m}$ solution of fractional differential equation
\begin{equation}\label{Eq:newproblem4}
\left\{
\begin{array}{ll}
\big(\fdabc{0}{\alpha} (y_{m}(t),\vartheta) \big)_{L^{2}(\Omega)} + a\big(y_{m}(t),~ \vartheta\big)_{L^{2}(\Omega)}  = \big(f(t),~\vartheta \big)_{L^{2}(\Omega)} & \text{for all} ~~\vartheta \in  \mathcal{V}_{m},\\[3mm]
y_{m}(x,0) = y_{m}^{0} &\text{for } x \in \Omega.
\end{array}
\right.
\end{equation}

\noindent Since $y_{m}(t) \in \mathcal{V}_{m}$, we have
\[
y_{m}(t) = \sum_{i=1}^{m}\big(y(t),w_{i}\big)_{L^{2}(\Omega)}w_{i} = \sum_{i=1}^{m}y_{i}(t)w_{i},
\]
which could still be written in an explicit form as
\begin{equation}\label{Eq:solution}
\begin{aligned}
y_{m}(x,t)  = \sum_{i=1}^{m}\bigg\{ \zeta_{i} E_{\alpha}\big[-\gamma_{i} t^{\alpha}\big]y_{i}^{0}
+ \frac{(1-\alpha) \zeta_{i}}{B(\alpha)}f_{i}(t) + K_{i} \int_{0}^{t}(t-s)^{\alpha-1} E_{\alpha,\alpha}\big[-\gamma_{i}(t-s)^{\alpha}\big]f_{i}(s)ds\bigg\}w_{i}.
\end{aligned}
\end{equation}

\noindent The next step is to show  $y_{m}(t)$ is a Cauchy sequence in the space $L^{2}\big((0,T);H_{0}^{1}(\Omega)\big)$. Let $m$ and $r$ be two natural numbers such that $r>m$. Then
\[
y_{r}(t) - y_{m}(t) = \sum_{i=m+1}^{r}y_{i}(t)w_{i}.
\]
\noindent But we also have that
\begin{equation}\label{long_computation_sequences1}
\begin{aligned}
& a\big(y_{r}(t)-y_{m}(t), y_{r}(t)-y_{m}(t) \big) = \sum_{i=m+1}^{r}\lambda_{i} \big(y_{i}(t)\big)^{2}\\
& = \sum_{i=m+1}^{r}\lambda_{i} \bigg\{ \zeta_{i} E_{\alpha}\big[-\gamma_{i} t^{\alpha}\big]y_{i}^{0}
+ \frac{(1-\alpha) \zeta_{i}}{B(\alpha)}f_{i}(t) + K_{i} \int_{0}^{t}(t-s)^{\alpha-1} E_{\alpha,\alpha}\big[-\gamma_{i}(t-s)^{\alpha}\big]f_{i}(s)ds\bigg\}^{2}\\
& \leq 3\sum_{i=m+1}^{r}\lambda_{i} \bigg[ \zeta_{i} E_{\alpha}\big[-\gamma_{i} t^{\alpha}\big]y_{i}^{0}\bigg]^{2} + 3\sum_{i=m+1}^{r}\lambda_{i} \bigg[\frac{(1-\alpha) \zeta_{i}}{B(\alpha)}f_{i}(t) \bigg]^{2} \\
&+ 3\sum_{i=m+1}^{r}\lambda_{i}K^{2}_{i} \bigg[ \int_{0}^{t}(t-s)^{\alpha-1} E_{\alpha,\alpha}\big[-\gamma_{i}(t-s)^{\alpha}\big]f_{i}(s)ds \bigg]^{2}.
\end{aligned}
\end{equation}
Hence
\[
\int_{0}^{T}\| y_{r}(t) - y_{m}(t) \|^{2}_{L^{2}\big((0,T);H_{0}^{1}(\Omega)\big)}dt  = \int_{0}^{T} a\big(y_{r}(t)-y_{m}(t), y_{r}(t)-y_{m}(t) \big) dt \leq 3(Z_{1} + Z_{2} + Z_{3}),
\]
where $Z_{1}$, $Z_{2}$ and $Z_{3}$ are defined respectively as
\[
\begin{aligned}
 & Z_{1}:= \sum_{i=m+1}^{r}\lambda_{i} \int_{0}^{T} \bigg[ \zeta_{i} E_{\alpha}\big[-\gamma_{i} t^{\alpha}\big]y_{i}^{0}\bigg]^{2}dt, \quad  Z_{2}:= \sum_{i=m+1}^{r}\lambda_{i} \int_{0}^{T} \bigg[\frac{(1-\alpha) \zeta_{i}}{B(\alpha)}f_{i}(t) \bigg]^{2} dt,\\
&Z_{3}:= \sum_{i=m+1}^{r}\lambda_{i}K^{2}_{i} \int_{0}^{T} \bigg[\int_{0}^{t}(t-s)^{\alpha-1} E_{\alpha,\alpha}\big[-\gamma_{i}(t-s)^{\alpha}\big]f_{i}(s)ds \bigg]^{2} dt.
\end{aligned}
\]

\noindent Now the goal is to estimate each term $Z_{1}$, $ Z_{2}$ and $Z_{3}$.\medskip

\noindent $\bullet$ {\sc Estimate of $ Z_{1}$.}
Using the estimate of the Mittag-Leffler function provided in Lemma~\ref{lem:mlfbdd} we have that
\begin{equation}\label{long_computation_Z1}
\begin{aligned}
Z_{1}  & = \sum_{i=1+m}^{r}\lambda_{i} \int_{0}^{T}  \bigg[ \zeta_{i} E_{\alpha}\big[-\gamma_{i} t^{\alpha}\big]u_{i}^{0}\bigg]^{2} dt \leq 2C^{2} \sum_{i=m+1}^{r}\lambda_{i} \zeta_{i}^{2} |u_{i}^{0}|^{2} \int_{0}^{T} \bigg(\frac{1}{1+\gamma_{i} t^{\alpha}} \bigg)^{2} dt \\
& \leq 2C^{2} T \sum_{i=m+1}^{r}\lambda_{i} \zeta_{i}^{2} |u_{i}^{0}|^{2}
\end{aligned}
\end{equation}
Now we find the estimate $ \lambda_{i}\zeta_{i}^{2} $ as
\[
\lambda_{i}\zeta_{i}^{2} ~ = \frac{\lambda_{i} B^{2}(\alpha)}{\big( B(\alpha) + (1-\alpha)\lambda_{i} \big)^{2}} \leq  ~~ \frac{\lambda_{i} B^{2}(\alpha)}{ (1-\alpha)^{2}\lambda_{i}^{2}} ~\leq  ~~ \frac{B^{2}(\alpha)}{(1-\alpha)^{2}\lambda_{i} } ~\leq  ~~ \frac{B^{2}(\alpha)}{(1-\alpha)^{2}\lambda_{1} },
\]
since $\lambda_{i}$ are increasing. We finally get the estimate of $Z_{1}$ as
\[
Z_{1}  \leq \frac{2C^{2} B^{2}(\alpha)}{\lambda_{1}(1-\alpha)^{2}}T\sum_{i=m+1}^{r} |u_{i}^{0}|^{2}.
\]

\noindent $\bullet$ {\sc Estimate of $ Z_{2}$.}

\begin{equation}\label{long_computation_Z2}
\begin{aligned}
Z_{2} & = 2\sum_{i=m+1}^{r}\lambda_{i} \int_{0}^{T} \bigg[\frac{(1-\alpha) \zeta_{i}}{B(\alpha)}f_{i}(t) \bigg]^{2} = \frac{2(1-\alpha)^{2}}{B(\alpha)^{2}}\sum_{i=m+1}^{r}\lambda_{i} \zeta_{i}^{2} \bigg(\int_{0}^{T} \big\|f_{i}(t) \big\|^{2}dt \bigg) \\
& \leq \frac{2B^{2}(\alpha)}{(1-\alpha)^{2}\lambda_{1} }\sum_{i=m+1}^{r}\bigg(\int_{0}^{T} \big|f_{i}(t) \big|^{2}dt \bigg) .
\end{aligned}
\end{equation}

\noindent $\bullet$ {\sc Estimate of $ Z_{3}$.}
By applying the Cauchy-Schwarz inequality and by setting $\delta = t-s$ we have that
\[
\begin{aligned}
Z_{3} & = 2\sum_{i=m+1}^{r}\lambda_{i}K_{i}^{2}\int_{0}^{T} \bigg[\int_{0}^{t}(t-s)^{\alpha-1} E_{\alpha,\alpha}\big[-\gamma_{i}(t-s)^{\alpha-1}\big]f_{i}(s)ds \bigg]^{2} dt\\
& \leq 2\sum_{i=m+1}^{r}\lambda_{i}K_{i}^{2} \int_{0}^{T} \left(\int_{0}^{t}\left( \delta^{\alpha-1} E_{\alpha,\alpha}\left[-\gamma_{i}\delta^{\alpha-1} \right]\right)^{2} d\delta \right) \left(\int_{0}^{t} \|f_{i}(s) \|^{2}ds\right) dt
\end{aligned}
\]
Now from Lemma~\ref{generalized} we have that
\[
\begin{aligned}
&\frac{d}{d\delta}  E_{\alpha,1}\big[-\gamma_{i} \delta^{\alpha}\big]  = \frac{d}{d\delta}  \mathcal{E}^{1}_{\alpha,1}\big[-\gamma_{i} \delta^{\alpha}\big] = \mathcal{E}^{2}_{\alpha,\alpha+1}\big[-\gamma_{i} \delta^{\alpha}\big] \\
&= - \gamma_{i} \delta^{\alpha-1} E_{\alpha,\alpha}\big[-\gamma_{i} \delta^{\alpha}\big] = \frac{1}{\alpha}   E_{\alpha,\alpha}\big[-\gamma_{i} \delta^{\alpha}\big] d \delta.
\end{aligned}
\]
Furthermore the estimate of $\lambda_{i}K_{i}^{2}$ is giving by
\[
\begin{aligned}
\lambda_{i} K_{i}^{2} &= \lambda_{i} \left(\frac{\alpha \zeta_{i}}{B(\alpha)\Gamma(\alpha)} + \frac{(1-\alpha) \gamma_{i}\zeta_{i}}{B(\alpha)}\right)^{2} \leq 2\left(\frac{\alpha^{2}}{B^{2}(\alpha)\Gamma^{2}(\alpha)}\lambda_{i} \zeta_{i}^{2} + \frac{(1-\alpha)^{2} }{B^{2}(\alpha)}\gamma_{i}^{2}\lambda_{i}\zeta_{i}^{2}\right)\\
& \leq \frac{2\alpha^{2}}{\lambda_{1}(1-\alpha)^{2} } \left(1+\frac{1}{\Gamma^{2}(\alpha)} \right)
\end{aligned}
\]
Then we get that
\begin{equation}\label{long_computation_Z3}
\begin{aligned}
Z_{3} & \leq \frac{4\alpha^{2}}{\lambda_{1}(1-\alpha)^{2} } \left( 1+ \frac{1}{\Gamma^{2}(\alpha)} \right) \sum_{i=m+1}^{r} \frac{1}{\gamma_{i} ^{2}}\int_{0}^{T} \left(\int_{0}^{t}\left( E_{\alpha,\alpha}\left[-\gamma_{i}\delta^{\alpha-1} \right]\right)^{2} d\delta \right) \left(\int_{0}^{t} \|f_{i}(s) \|^{2}ds\right) dt\\
& \leq \frac{4C^{2}}{\lambda_{1}} \left( 1+ \frac{1}{\Gamma^{2}(\alpha)} \right)T \sum_{i=m+1}^{r} \int_{0}^{T} \left(\int_{0}^{T} \|f_{i}(t) \|^{2}dt\right) dt\\
& \leq \frac{4C^{2}T^{2}}{\Gamma^{2}(\alpha)} \left(\frac{\Gamma^{2}(\alpha)+ 1}{\lambda_{1}} \right) \sum_{i=m+1}^{r} \int_{0}^{T} \|f_{i}(t) \|^{2} dt \\
\end{aligned}
\end{equation}

\noindent Finally, adding together the estimates $Z_{1}$, $Z_{2}$, and $Z_{3}$  we get
\[
\begin{aligned}
& \int_{0}^{T}\big\| y_{r}(x,t) - y_{m}(x,t) \big\|^{2}_{L^{2}\big((0,T);H_{0}^{1}(\Omega)\big)} dt
\leq  \frac{6C^{2} B^{2}(\alpha)}{\lambda_{1}(1-\alpha)^{2}}T\sum_{i=m+1}^{r} |y_{i}^{0}|^{2} \\
&+ \frac{6B^{2}(\alpha)}{(1-\alpha)^{2}\lambda_{1} }\sum_{i=m+1}^{r}\bigg(\int_{0}^{T} \big|f_{i}(t) \big|^{2}dt \bigg) + \frac{12C^{2}T^{2}}{\Gamma^{2}(\alpha)} \left(\frac{\Gamma^{2}(\alpha)+ 1}{\lambda_{1}} \right)\sum_{i=m+1}^{r} \int_{0}^{T} \|f_{i}(t) \|^{2} dt
\end{aligned}
\]
Thus  we obtain
\[
\begin{aligned}
\int_{0}^{T}\big\| y_{r} - y_{m}\big\|_{H_{0}^{1}(\Omega)} &\leq \mathcal{C}_{1}(\alpha,\lambda_{1},T)\left(\sum_{i=m+1}^{r} |y_{i}^{0}|^{2} \right)^{1/2} \\
& + \mathcal{C}_{2}(\alpha,\lambda_{1},T) \left( \sum_{i=m+1}^{r}  \int_{0}^{T} \big| f_{i}(t)\big|^{2}dt \right)^{1/2},
\end{aligned}
\]
with
\[
\mathcal{C}_{1}(\alpha,\lambda_{1},T) =  \frac{C B(\alpha)}{(1-\alpha)}\sqrt{\frac{6T}{\lambda_{1}}}, \qquad \mathcal{C}_{2}(\alpha,\lambda_{1},T) = \sqrt{\frac{6B^{2}(\alpha)}{(1-\alpha)^{2}\lambda_{1} } + \frac{12C^{2}T^{2}}{\Gamma^{2}(\alpha)} \left(\frac{\Gamma^{2}(\alpha)+ 1}{\lambda_{1}} \right)}.
\]

\noindent As in the previous case, the norm in $L^{2}(\Omega)$ is computed as

\begin{equation}\label{ll}
\begin{aligned}
&\int_{0}^{T}\big\| y_{r} - y_{m}\big\|_{L^{2}(\Omega)}^{2}dt = \int_{0}^{T} \sum_{i=m+1}^{r} \big(y_{i}(t)\big)^{2} dt\\
& = \int_{0}^{T}\sum_{i=m+1}^{r} \bigg\{ \zeta_{i} E_{\alpha}\big[-\gamma_{i} t^{\alpha}\big]y_{i}^{0}
+ \frac{(1-\alpha) \zeta_{i}}{B(\alpha)}f_{i}(t) + K_{i} \int_{0}^{t}(t-s)^{\alpha-1} E_{\alpha,\alpha}\big[-\gamma_{i}(t-s)^{\alpha}\big]f_{i}(s)ds\bigg\}^{2} dt\\
& \leq 3 \int_{0}^{T}\sum_{i=m+1}^{r}\bigg[ \zeta_{i} E_{\alpha}\big[-\gamma_{i} t^{\alpha}\big]y_{i}^{0}\bigg]^{2}dt + 3\int_{0}^{T} \sum_{i=m+1}^{r} \bigg[\frac{(1-\alpha) \zeta_{i}}{B(\alpha)}f_{i}(t) \bigg]^{2}dt \\
&+ 3\int_{0}^{T}\sum_{i=m+1}^{r}K^{2}_{i} \bigg[ \int_{0}^{t}(t-s)^{\alpha-1} E_{\alpha,\alpha}\big[-\gamma_{i}(t-s)^{\alpha}\big]f_{i}(s)ds \bigg]^{2}dt.
\end{aligned}
\end{equation}
\noindent From the previous estimations in the first step of the proof, we have that
\begin{equation}\label{ll1}
\begin{aligned}
\int_{0}^{T}\big\| y_{r} - y_{m}\big\|_{L^{2}(\Omega)}^{2}dt &=
\frac{6C^{2} B^{2}(\alpha)}{(1-\alpha)^{2}\lambda_{1}^{2}}T \sum_{i=m+1}^{r} \|y_{i}^{0}\|^{2}
+ \frac{6}{\lambda_{1}^{2}} \sum_{i=m+1}^{r} \int_{0}^{T} \big\|f_{i}(t) \big \|^{2} \\
&+ \frac{12C^{2}T^{2}}{\Gamma^{2}(\alpha)} \left(\frac{\Gamma^{2}(\alpha)+ 1}{\lambda_{1}^{2}} \right) \sum_{i=m+1}^{r} \int_{0}^{T} \|f_{i}(t) \|^{2} dt
\end{aligned}
\end{equation}

\noindent Thus we obtain
\[
\begin{aligned}
&\sup_{t\in[0,T]}\big\| y_{r}- y_{m} \big\|_{L^{2}(\Omega)} \\
& \leq \mathcal{C}_{3}(\alpha,\lambda_{1},T)\bigg(\sum_{i=m+1}^{r} |y_{i}^{0}|^{2} \bigg)^{1/2} + ~ \mathcal{C}_{4}(\alpha,\lambda_{1},T) \left( \sum_{i=m+1}^{r}  \int_{0}^{T} \big| f_{i}(t)\big|^{2}dt \right)^{1/2},
\end{aligned}
\]
with
\[
\mathcal{C}_{3}(\alpha,\lambda_{1},T) =  \frac{C B(\alpha)}{\lambda_{1}(1-\alpha)}\sqrt{6T}, \qquad \mathcal{C}_{4}(\alpha,\lambda_{1},T) = \sqrt{\frac{6}{\lambda_{1}^{2}} + \frac{12C^{2}T^{2}}{\Gamma^{2}(\alpha)} \left(\frac{\Gamma^{2}(\alpha)+ 1}{\lambda_{1}^{2}} \right)}.
\]

\noindent From the fact that $ y^{0} \in H_{0}^{1}(\Omega) $ and $f \in  L^{2}(\mathcal{Q})$, we have that
\[
\lim_{m,r\to +\infty} \left(\sum_{i=m+1}^{r} |y_{i}^{0}|^{2} \right)^{1/2} = 0, \qquad
\lim_{m,r\to +\infty} \left( \sum_{i=m+1}^{r}  \int_{0}^{T} \big\| f_{i}(t)\big\|^{2}dt \right)^{1/2} = 0.
\]
Hence we have that
\[
\lim_{m,r\to +\infty} \int_{0}^{T}\big\| y_{r} - y_{m} \big\|_{L^{2}\left((0,T);H_{0}^{1}(\Omega) \right)}^{2} = 0,
\]
and
\[
\lim_{m,r\to +\infty} \int_{0}^{T}\big\| y_{r}- y_{m}\big\|_{L^{2}(\Omega)}^{2} = 0,
\]
which means that the sequence $(y_{m})$ is  Cauchy sequences in  $L^{2}\left((0,T);H_{0}^{1}(\Omega) \right)$ and  $\mathcal{C}\left([0,T];L^{2}(\Omega) \right)$. \medskip

Therefore we have that 
\[
\begin{array}{lllllll}
y_{m}&\rightharpoonup&y&\hbox{ in }&L^{2}\left((0,T);H_{0}^{1}(\Omega) \right),\\
y_{m}&\rightharpoonup&y&\hbox{ in }&\mathcal{C}\left([0,T];L^{2}(\Omega) \right).
\end{array}
\]
Now proceeding as in \cite{BaleanuJosephMophou}, one can prove by interpretation that the solution $y$ is equivalent to the initial problem \eqref{Eq:problem}.
\noindent Furthermore, setting $\Lambda_{1} = \max\left(\mathcal{C}_{1}(\alpha,\lambda_{1},T),~ \mathcal{C}_{2}(\alpha,\lambda_{1},T)\right)$ and  $\Lambda_{2} = \max\left(\mathcal{C}_{3}(\alpha,\lambda_{1},T),~ \mathcal{C}_{4}(\alpha,\lambda_{1},T)\right)$, an using the estimates obtained at the beginning of the proof,  we have that
\[
\begin{array}{lll}
\|y\|_{L^2((0,T);H^1_0(\Omega))}&=&\displaystyle\left(\int_{0}^{T}\big\| y\|_{H^1_0(\Omega)} dt\right)^{1/2}\\
&\leq& \mathcal{C}_{1}(\alpha,\lambda_{1},T)\|y^{0}\|_{L^2(\Omega)}  +
\mathcal{C}_{2}(\alpha,\lambda_{1},T) \| f\|_{L^2(Q)}\\
&\leq& \Lambda_{1} \left( \|y^{0}\|_{H_{0}^{1}(\Omega)\big)} + \|f\|_{L^{2}(\mathcal{Q})}\right),
\end{array}
\]
and 
\[
\begin{array}{lll}
\sup_{t\in[0,T]}\big\| y \big\|_{L^{2}(\Omega)} & \leq& \mathcal{C}_{3}(\alpha,\lambda_{1},T) \|y\|_{L^2(\Omega)}+ \mathcal{C}_{4}(\alpha,\lambda_{1},T) \|f\|_{L^2(Q)}\\
&\leq& \Lambda_{2} \left( \|y^{0}\|_{L^{2}(\Omega)} + \|f\|_{L^{2}(Q)}\right)
\end{array}
\]
Hence we get the desired result, which ends the proof.\medskip
\end{proof}
\begin{rem}\label{existence_range}
We would like to draw the attention of the reader that in the case for the Caputo fractional-time derivative, the existence and uniqueness of solutions for the fractional diffusion equation of type \eqref{Eq:problem}, were obtained under the conditions that $\alpha \in (1/2,1)$ \cite{Yamamoto,Mophou,BaleanuJosephMophou}. But here in our case using the Atangana-Baleanu fractional-time derivative, we got the existence and uniqueness of solutions for the fractional diffusion equation of type \eqref{Eq:problem} for all $\alpha \in (0,1)$.
\end{rem}
\noindent Next we state two corollaries which give the estimates of the solution $y(t)$ in $L^{2}\left((0,T),H^{2}(\Omega)\cap H^{1}_{0}(\Omega)\right)$ for the case $f = 0$ and in the case where the initial data $y^{0} = 0$ respectively in the problem \eqref{Eq:problem} or \eqref{Eq:newproblem1}.\medskip

\noindent $\bullet$ {\sc First case: $f = 0$}.

\noindent For $f = 0$ in the problem \eqref{Eq:problem} or \eqref{Eq:newproblem1} we have the following Corollary.
\begin{cor}\label{coro1}
\noindent Let $\alpha \in (0,1)$. Assume $f =0$, $y^{0} \in L^{2}(\Omega)$ and the bilinear form $a(\cdot,\cdot)$ defined as in \eqref{bilinear_form}. Then the problem \eqref{Eq:problem} or \eqref{Eq:newproblem1} has a unique solution $y \in L^{2}\left((0,T),H^{2}(\Omega)\cap H^{1}_{0}(\Omega)\right)$ given by
\[
\begin{aligned}
y(x,t)  = \sum_{i=1}^{+\infty}\bigg\{ \zeta_{i} E_{\alpha}\big[-\gamma_{i} t^{\alpha}\big]y_{i}^{0}
\bigg\}w_{i},
\end{aligned}
\]
where the constant $\zeta_{i}$ is defined as in \eqref{contants}. Moreover $y$ satisfies the bounded conditions
\begin{equation}\label{aaai}
\int_{0}^{T}\big \| y \big\|_{H^{2}(\Omega)}^{2}dt = \|y\|^{2}_{L^{2}\left((0,T),H^{2}(\Omega)\right)} \leq \Lambda_{3} \|y^{0}\|_{H_{0}^{2}(\Omega)},
\end{equation}
where $\Lambda_{3} = \frac{CB(\alpha)}{(1-\alpha)}\sqrt{T}$.
\end{cor}

\begin{proof}

\noindent Following the same procedure as in the first step for the proof of the Theorem~\ref{existence_uniqueness}, but now in $H_{0}^{2}(\Omega)$, the solution of the problem \eqref{Eq:problem} or \eqref{Eq:newproblem1} for $f = 0$ is given by
\[
\begin{aligned}
y(t)  = \sum_{i=1}^{+\infty}\bigg\{ \zeta_{i} E_{\alpha}\big[-\gamma_{i} t^{\alpha}\big]y_{i}^{0}
\bigg\}w_{i},
\end{aligned}
\]
where the constant $\zeta_{i}$ is defined as in \eqref{contants}. Furthermore the estimate of $y$ can be obtained as follows
\begin{equation}\label{llll}
\begin{aligned}
\int_{0}^{T} \big \| y \big\|^{2}_{H^{2}(\Omega)} dt &= \int_{0}^{T}  \sum_{i=m+1}^{r}\lambda_{i}^{2} \big(y_{i}(t)\big)^{2}  dt = \int_{0}^{T} \sum_{i=m+1}^{r}\lambda_{i}^{2} \bigg\{ \zeta_{i} E_{\alpha}\big[-\gamma_{i} t^{\alpha}\big]y_{i}^{0} \bigg\}^{2} dt\\
& \leq \sum_{i=m+1}^{r}\lambda_{i}^{2}\zeta_{i}^{2} \int_{0}^{T} \left( E^{2}_{\alpha}\big[-\gamma_{i} t^{\alpha}\big]\right)^{2} \big\|y_{i}^{0} \big \|^{2} dt\\
& \leq \frac{C^{2}B^{2}(\alpha)}{(1-\alpha)^{2}}T\sum_{i=m+1}^{r}\big\|y_{i}^{0} \big \|^{2} \\
& \leq \mathcal{C}_{5}(\alpha,T) \sum_{i=m+1}^{r}\big\|y_{i}^{0} \big \|^{2}.
\end{aligned}
\end{equation}
\end{proof}

\noindent $\bullet$ {\sc Second case: $y^{0}_{i} = 0$}. \medskip

\noindent For $y^{0}_{i} = 0$ in the problem \eqref{Eq:problem} or \eqref{Eq:newproblem1} we have the following Corollary.
\begin{cor}\label{coro2}
\noindent Let $\alpha \in (0,1)$. Assume $y^{0}_{i} = 0$, and $f \in L^{2}(Q)$. Then the problem \eqref{Eq:problem} or \eqref{Eq:newproblem1} has a unique solution $y \in L^{2}\left((0,T),H^{2}(\Omega)\cap H^{1}_{0}(\Omega)\right)$ given by
\[
\begin{aligned}
y(x,t)  = \sum_{i=1}^{+\infty}\bigg\{\frac{(1-\alpha) \zeta_{i}}{B(\alpha)}f_{i}(x, t) + K_{i} \int_{0}^{t}(t-s)^{\alpha-1} E_{\alpha,\alpha}\big[-\gamma_{i}(t-s)^{\alpha}\big]f_{i}(s)ds\bigg\}w_{i},
\end{aligned}
\]
where the constants $\gamma_{i}$, $\zeta_{i}$ and $K_{i}$ are defined as in \eqref{contants}. Moreover $y$ satisfies the bounded conditions
\begin{equation}\label{aaaii}
\big \| y \big\|_{L^{2}\left((0,T),H^{2}(\Omega)\right)} \leq T\Lambda_{4} \|f\|_{L^{2}(Q)},
\end{equation}
where $\Lambda_{4} = \sqrt{\left(2 + 4 C^{2}T^{2}\left[1 + \frac{1}{\Gamma^{2}(\alpha)} \right] \right)}$,
\end{cor}

\begin{proof}
Let $\alpha \in (0,1)$. Assume $y^{0}_{i} = 0$, and $f \in L^{2}(Q)$. Then the problem \eqref{Eq:problem} or \eqref{Eq:newproblem1} has a unique solution $y \in H^{2}(\Omega)\big)$ given by
\[
\begin{aligned}
y(x,t)  = \sum_{i=1}^{+\infty}\bigg\{\frac{(1-\alpha) \zeta_{i}}{B(\alpha)}f_{i}(t) + K_{i} \int_{0}^{t}(t-s)^{\alpha-1} E_{\alpha,\alpha}\big[-\gamma_{i}(t-s)^{\alpha}\big]f_{i}(s)ds\bigg\}w_{i},
\end{aligned}
\]
where the constants $\gamma_{i}$, $\zeta_{i}$ and $K_{i}$ are defined as in \eqref{contants}. Furthermore we compute the estimate of $y$ in the space $H^{2}(\Omega)$ as
\begin{equation}\label{lllli}
\begin{aligned}
\int_{0}^{T} \big \| y \big\|^{2}_{H^{2}(\Omega)} dt &= \int_{0}^{T}  \sum_{i=m+1}^{r}\lambda_{i}^{2} \bigg\{\frac{(1-\alpha) \zeta_{i}}{B(\alpha)}f_{i}(t) + K_{i} \int_{0}^{t}(t-s)^{\alpha-1} E_{\alpha,\alpha}\big[-\gamma_{i}(t-s)^{\alpha}\big]f_{i}(s)ds\bigg\}^{2}  dt  \\
&\leq 2\sum_{i=m+1}^{r}\lambda_{i}^{2} \int_{0}^{T} \bigg[\frac{(1-\alpha) \zeta_{i}}{B(\alpha)}f_{i}(t) \bigg]^{2}\\
& + 2\sum_{i=m+1}^{r}\lambda_{i}^{2}K_{i}^{2} \int_{0}^{T} \left(\int_{0}^{t}\left((t-s)^{\alpha-1} E_{\alpha,\alpha}\left[-\gamma_{i}(t-s)^{\alpha-1} \right]\right)^{2} ds \right) \left(\int_{0}^{t} \|f_{i}(s) \|^{2}ds\right) dt \\
& \leq 2\sum_{i=m+1}^{r} \int_{0}^{T} \big\|f_{i}(t) \big\|^{2} + 4C^{2}T\left(1 + \frac{1}{\Gamma^{2}(\alpha)} \right) \sum_{i=m+1}^{r} \int_{0}^{T} \big\|f_{i}(t) \big\|^{2} dt \\
& \leq \mathcal{C}_{6}\sum_{i=m+1}^{r} \int_{0}^{T} \big\|f_{i}(t) \big\|^{2} dt
\end{aligned}
\end{equation}
\end{proof}

\noindent Now, we consider the following fractional differential equation:
\begin{equation}\label{Eq:fde_adjoint}
\left\{
\begin{array}{ll}
-\fdabc{0}{\alpha} \eta(t) - \mathcal{A}\eta(t) = g(t), \qquad t\in (0,T),\\[3mm]
\eta(T) = 0,
\end{array}
\right.
\end{equation}
where $\alpha \in (0,1)$, $g \in L^{2}(Q)$.
We state the following Lemma.
\begin{lem}\label{adjoin}
Let $\alpha \in (0,1)$, $g \in L^{2}(Q)$. Then problem \eqref{Eq:fde_adjoint} has a unique solution $\eta \in L^{2}(Q)$ given by
\begin{equation}\label{Eq:solution_adjoin}
\begin{aligned}
\eta_{i}(x,t)  = \sum_{i=1}^{+\infty} \left\lbrace \frac{(1-\alpha) \zeta_{i}}{B(\alpha)}g(t) + K_{i} \int_{0}^{t}(t-s)^{\alpha} E_{\alpha,\alpha}\big[-\gamma_{i}(t-s)^{\alpha}\big]g(s)ds \right\rbrace w_{i}.
\end{aligned}
\end{equation}
with the same constants as defined in \eqref{contants}. Moreover

\begin{equation}\label{bbb}
\big\|\eta_{i} \big\|_{L^{2}\big((0,T);H_{0}^{1}(\Omega)\big)} \leq \sqrt{\mathcal{C}_{2}}  \|g\|_{L^{2}(Q)}.
\end{equation}
\end{lem}

\begin{proof}
\noindent Using Lemma~\ref{shift}, we found that the backwards in time with the Atangana-Baleanu fractional-time derivative at the based point $T$ is equivalently written as a forward in time operator with the Atangana-Baleanu fractional-time derivative at a based point zero time a negative sign. More precisely,
\[
\fdabc{T}{\alpha}~\widetilde{\eta}(t) = -\fdabc{0}{\alpha}~\widetilde{\eta}(t).
\]
\noindent Using this relation and by making the change of variable $t \to T-t$ in \eqref{Eq:fde_adjoint}, we obtain the fractional differential equation
\[
\left\{
\begin{array}{ll}
\fdabc{0}{\alpha} \widetilde{\eta}(t) - \mathcal{A}\widetilde{\eta}(t) = \widetilde{g}(t), \qquad T-t\in [0,T],\\[3mm]
\eta(0) = 0,
\end{array}
\right.
\]
with $\widetilde{g}(t) = g(T-t)$. This fractional differential equation can still be written in the short form as
\begin{equation}\label{Eq:fde_adjoint_a}
\left\{
\begin{array}{ll}
\fdabc{0}{\alpha} \eta(\tau) - \mathcal{A}\eta(\tau) = g(\tau), \qquad \tau \in [0,T],\\[3mm]
\eta(0) = 0,
\end{array}
\right.
\end{equation}
\noindent Next to show that \eqref{Eq:solution_adjoin} holds, we just apply the Laplace transform on \eqref{Eq:fde_adjoint_a} and get the inverse following the steps as in \cite[Proposition $3.3$]{DjidaIvan_Abdon_regularity}, and by taking into account the fact that $\eta(0) = 0$. Hence one get
\[
\begin{aligned}
\eta(t)  =  \frac{(1-\alpha) \zeta_{i}}{B(\alpha)}g(t) + K_{i} \int_{0}^{t}(t-s)^{\alpha} E_{\alpha,\alpha}\big[-\gamma_{i}(t-s)^{\alpha}\big]g(s)ds,
\end{aligned}
\]
where $\zeta_{i}$, $\gamma_{i}$ and $\mathcal{C}_{2}$ were defined respectively in \eqref{contants} and in  \eqref{ccc} \medskip

\noindent Finally proceeding as for the proof of \eqref{aaa} in Theorem \ref{existence_uniqueness} once can obtain
\[
\big\|\eta_{i} \big\|_{L^{2}\big((0,T);H_{0}^{1}(\Omega)\big)} \leq \sqrt{\mathcal{C}_{2}}  \|g\|_{L^{2}(Q)}.
\]
\end{proof}

\section{Optimal control problem}\label{control}

\noindent In this section, we want to approach the state $y(v)$ of the problem \eqref{Eq:problem} by a desired state
$z_{d}$ in controlling $v$. \medskip

\noindent Let $v\in L^{2}(Q)$. As it has been shown in the existence and uniqueness part, the problem \eqref{Eq:problem} have a unique solution $y= y(v)$ which belongs to $L^{2}(Q)$. In this way, we can define can define the functional as
\begin{equation}\label{functional}
\mathcal{J}(v) = \big\|y(v)-z_{d}\big \|^{2}_{L^{2}(Q)} + \frac{\mathcal{N}}{2}\|v\|^{2}_{L^{2}(Q)},
\end{equation}
where $z_{d} \in L^{2}(Q)$ and $ \mathcal{N}>0$. We are interested on finding the control parameter, say, $\hat{u} \in L^{2}(Q)$ such that the functional
\begin{equation}\label{functional_minimizer}
\mathcal{J}(\hat{u}) = \inf_{v\in L^{2}(Q)} \mathcal{J}(v).
\end{equation}
\noindent Now we state the following Proposition.
\begin{prop}
Let consider the problem given by \eqref{Eq:problem}. Then there exists a unique optimal control $u$ such that \eqref{functional_minimizer} holds true.
\end{prop}

\begin{proof}
The functional $\mathcal{J}$  is continuous, coercive and convex. Therefore, there exist a control called $\hat{u}$. By calling $\hat{y} = \hat{y}(\hat{u})$, the associated state of $\hat{u}$ then its come that, $\hat{y}$ satisfies

\begin{equation}\label{Eq:associate_states}
\left\{
\begin{array}{ll}
\fdabc{0}{\alpha} \hat{y}- \Delta \hat{y} = \hat{u} & \text{for } (x,t)\in Q,\\
\hat{y}(x,0) = \hat{y}^{0}(x) &\text{for } x \in \Omega,\\
\hat{y}(x,t) = 0 & \text{on} ~\partial \Omega \times (0,T),
\end{array}
\right.
\end{equation}
where $y \in L^{2}((0,T); H^{1}_{0}(\Omega)) \cap \mathcal{C}([0,T]; L^{2}(\Omega))$. \medskip
\end{proof}

\begin{rem}
\noindent Another way of proving this result is to use minimizing sequences. The reader interested could follows the idea proposed by the second author in \cite{Mophou}.
\end{rem}

\noindent In the following we are concern with the interpretation of the Euler-Lagrange first order optimality condition with an adjoint problem defined by means of the backward Atangana-Baleanu fractional-time derivative. As a result, we obtained an optimality system for the optimal control given by the following Theorem.\medskip

\begin{sa}\label{optimality_system}
If $\hat{u}$ is a solution of \eqref{functional_minimizer}, then there exist $\eta$ such that 
\[
(\hat{u},\hat{y},\eta) \in L^{2}(Q) \times L^{2}\left((0,T); H^{1}_{0}(\Omega)\right) \times L^{2}\left((0,T); H^{2}(\Omega) \cap  H^{1}_{0}(\Omega)\right),
\]
satisfies the following optimality system
\begin{equation}\label{Eq:associate_statestheo}
\left\{
\begin{array}{ll}
-\fdabc{0}{\alpha} \hat{y}- \Delta \hat{y} = \hat{u} & \text{in} ~ Q,\\
\qquad \qquad \hat{y}(x,0) = \hat{y}^{0} &\text{in} ~ \Omega,\\
\qquad \qquad \hat{y}(x,t) = 0 & \text{on} ~\partial \Omega \times (0,T).
\end{array}
\right.
\end{equation}
\begin{equation}\label{Eq:associate_states_eta}
\left\{
\begin{array}{ll}
-\fdabc{T}{\alpha} \eta - \Delta \eta &= \hat{y}-z_{d} \qquad  \text{in}~ Q,\\
~~~\eta(T)& = 0 \qquad \qquad  \text{in } ~ \Omega,\\
~~~\eta &= 0  \qquad \qquad \text{on} ~\partial \Omega \times (0,T).
\end{array}
\right.
\end{equation}
\begin{equation}
\hat{u} = -\frac{\eta}{N} \quad \text{in}~~Q.
\end{equation}
\end{sa}

\begin{proof}
In order to prove the theorem, we start by expressing the Euler-Lagrange optimality conditions which characterize the optimal control $\hat{u}$. We means
\begin{equation}
\frac{d}{d\mu} \mathcal{J}\left(\hat{u} + \mu (v-\hat{u})\right) \vert_{\mu=0}, \quad \text{for all } \phi \in L^{2}(Q).
\end{equation}
The state $z\left(v-\hat{u}\right) \in L^{2}\left((0,T); H^{2}(\Omega) \cap  H^{1}_{0}(\Omega)\right) $  associated to the control $\left(v-\hat{u}\right) \in L^{2}(Q) $ is a  solution of
\begin{equation}\label{Eq:associate_statestheqbb}
\left\{
\begin{array}{ll}
\fdabc{0}{\alpha} z - \Delta z &= v-\hat{u} \quad  \text{in} ~~Q,\\
\qquad z\left(0,v-\hat{u}\right) &= 0 \qquad~~ \text{in} \Omega,\\
\qquad z &= 0 \qquad ~~\text{on} ~\partial \Omega \times (0,T).
\end{array}
\right.
\end{equation}
Now we find the following behaviour of our functional when $\mu$ converges to $0$.
\begin{equation}
\lim_{\mu \to 0} \frac{\mathcal{J}\left(\hat{u} + \mu (v-\hat{u})\right)-\mathcal{J}(\hat{u})}{\mu} = 0, \quad \text{for all}~ v \in L^{2}(Q).
\end{equation}
Using the functional defined in \eqref{functional} we have that
\[
\begin{aligned}
 \mathcal{J}\left(\hat{u} + \mu (v-\hat{u}) \right) & = \frac{\mathcal{N}}{2} (\hat{u} + \mu (v-\hat{u})(\hat{u} + \mu (v-\hat{u})\\
&+ \frac{1}{2} \left(\hat{y}(\hat{u}) - z_{d} + \mu\hat{y}(v-\hat{u}) \right) \left(\hat{y}(\hat{u}) - z_{d} + \mu\hat{y}(v-\hat{u}) \right)\\
&= \frac{\mathcal{N}}{2} \|\hat{u}\|^{2}_{L^{2}(Q)} + \mu \mathcal{N}\hat{u}(v-\hat{u}) + \mu^{2}\frac{\mathcal{N}}{2} \|v-\hat{u}\|^{2}_{L^{2}(Q)} \\
& + \frac{1}{2}\|\hat{y}(\hat{u}) - z_{d} )\|^{2}_{L^{2}(Q)} + \mu (\hat{y}(\hat{u}) - z_{d})\left(\hat{y}(v-\hat{u})\right) \\
& +  \frac{\mu^{2}}{2} \|\hat{y}(v-\hat{u})\|^{2}_{L^{2}(Q)}
\end{aligned}
\]
Hence we get that
\[
\begin{aligned}
\lim_{\mu \to 0} \frac{\mathcal{J}\left(\hat{u} + \mu (v-\hat{u}) \right)-\mathcal{J}(\hat{u})}{\mu} & =  z\left(\hat{y}(\hat{u}) - z_{d} \right)  + \mathcal{N}\hat{u}\left(v-\hat{u}\right) = 0, \quad \text{for all } v \in L^{2}(Q).
\end{aligned}
\]
Now integrating over $Q$ we get that
\begin{equation}\label{final_functional}
\int_{0}^{T}\int_{\Omega} z\left(\hat{y}(\hat{u}) - z_{d} \right) dx~dt + \mathcal{N} \int_{0}^{T}\int_{\Omega} \hat{u}\left(v-\hat{u}\right)  dx~dt = 0 \qquad \text{for all} ~~ v \in L^{2}(Q).
\end{equation}
In order to make the interpretation~\eqref{final_functional}, we consider the adjoint state equation
\begin{equation}\label{Eq:associate_states_eta_a}
\left\{
\begin{array}{ll}
-\fdabc{T}{\alpha} \eta - \Delta \eta &= \hat{y}-z_{d} \qquad  \text{in}~ Q,\\
~~~\eta(T)& = 0 \qquad \qquad  \text{in } ~ \Omega,\\
~~~\eta &= 0  \qquad \qquad \text{on} ~\partial \Omega \times (0,T).
\end{array}
\right.
\end{equation}
But since $\hat{y}-z_{d} \in L^{2}((0,T),H_{0}^{1}(\Omega))$, applying the Lemma~\ref{adjoin}, we deduce that the problem given by \eqref{Eq:associate_states_eta_a} has a unique solution in $L^{2}\left((0,T);H^{2}(\Omega) \cap H_{0}^{1}(\Omega)\right)$. In this scenario, by multiplying the state equation system \eqref{Eq:associate_statestheqbb} by the solution $\eta$ of \eqref{Eq:associate_states_eta_a}, we obtain by \ref{corweak} that
\[
\begin{aligned}
\int_{0}^{T}\int_{\Omega}\left(\fdabc{0}{\alpha} z - \Delta z \right) \eta dx~dt &=  \int_{0}^{T}\int_{\Omega} \left(-\fdabc{T}{\alpha} \eta - \Delta \eta \right)z dx~dt \\
 &= \int_{0}^{T}\int_{\Omega} \left( \hat{y}(u)- z_{d} \right) z dx~dt.
\end{aligned}
\]
Hence one get the following relation from \eqref{Eq:associate_statestheqbb} and \eqref{final_functional}, we deduce that
\[
\int_{0}^{T}\int_{\Omega}\left(v-\hat{u}\right)~\eta~dx~dt =  -\mathcal{N}\int_{0}^{T}\int_{\Omega} \left(v-\hat{u}\right) \hat{u}~dx~dt \qquad \text{for all} ~~v \in L^{2}(Q).
\]
We then get finally that
\[
\hat{u} = -\frac{\eta}{\mathcal{N}} \qquad \text{in}~~ Q.
\]
\end{proof}

\section{Conclusion and remarks}
In this work the existence and uniqueness of solution in $L^{2}(Q)$ were proved by means of a spectral argument and the control problem has also been studied. We have shown in the first step that the existence of solution were obtained for all values of the fractional parameter $\alpha \in (0,1)$, contrary to the Caputo and Riemann fractional-time derivative where the existence and uniqueness results were obtained for $\alpha \in (1/2,1)$ with $y^{0} \equiv 0$. This reveals one particularity of the fractional-time derivative with the nonsingular Mittag-Leffler function. Moreover, we have also 
shown that one can approach the state $y(v)$ of \eqref{Eq:problem} by a desired state $z_{d}$ by controlling $v$ and compute the control $\hat{u}$ using the algorithm given by the optimality system and following the method of Lagrange. Let us also mention that by proceeding as in this work, one can also study the same problem by considering the Atangana-Baleanu fractional-time derivative in the Riemann-Liouville sense.

\bibliographystyle{plain}

\end{document}